\documentstyle[12pt,amssymb,amscd]{amsart}

\newtheorem{thm}{Theorem}[section]
\newtheorem{prop}[thm]{Proposition}
\newtheorem{lma}[thm]{Lemma}

\newenvironment{defn}{\trivlist
\item[{\bf\hspace{1.88pt}\addtocounter{thm}{1}
Definition~\arabic{section}.\arabic{thm}.}]}{\endtrivlist}

\def\indlimit{{\displaystyle \lim_{\longrightarrow} }}

\def\resultspace{{\vskip 0.075in}}

\def\chisub#1 {{\chi_{\lower2.5pt\hbox{$\scriptstyle #1$}}}}

\def\symbup#1{\Big\uparrow\rlap{$\vcenter{\hbox{$\scriptstyle #1$}}$}}

\def\symbse#1{\!\!\!\!\searrow\rlap{$\vcenter{\hbox{$\scriptstyle #1$}}$}}

\def\lsymbse#1{\llap{$\vcenter{\hbox{$\scriptstyle #1$}}$}\!\!\searrow }

\def\ZZ{\Bbb{Z}}
\def\IC{\Bbb{C}}
\def\IQ{\Bbb{Q}}

\def\angb#1{\langle #1 \rangle }

\numberwithin{equation}{section}
\begin{document}

\bibliographystyle{amsalpha209}

\newcommand{\ca}{{\cal A}}
\newcommand{\cb}{{\cal B}}
\newcommand{\cc}{{\cal C}}
\newcommand{\cd}{{\cal D}}
\newcommand{\cf}{{\cal F}}
\newcommand{\cg}{{\cal G}}
\newcommand{\ch}{{\cal H}}
\newcommand{\ck}{{\cal K}}
\newcommand{\cl}{{\cal L}}
\newcommand{\cn}{{\cal N}}
\newcommand{\cm}{{\cal M}}
\newcommand{\cs}{{\cal S}}
\newcommand{\ct}{{\cal T}}
 
\newcommand{\ot}{\otimes} 
\newcommand{\op}{\oplus}  
\newcommand{\pr}{\prime}

\newcommand{\matrrccc}[6]{\mbox{$ \left[ \begin{array}{ccc}
	    #1&#2&#3\\ #4&#5&#6 \end{array} \right] $}}
\newcommand{\matrccc}[3]{\mbox{$ \left[ \begin{array}{ccc}
	    #1&#2&#3 \end{array} \right] $}}
\newcommand{\matrcccc}[4]{\mbox{$ \left[ \begin{array}{cccc}
	    #1&#2&#3&#4 \end{array} \right] $}}
\newcommand{\matrrrc}[3]{\mbox{$ \left[ \begin{array}{c}
	    #1 \\ #2 \\ #3 \end{array} \right] $}}
\newcommand{\matrrrrc}[4]{\mbox{$ \left[ \begin{array}{c}
	    #1 \\ #2 \\ #3 \\ #4 \end{array} \right] $}}
\newcommand{\matrrrcc}[6]{\mbox{$ \left[ \begin{array}{cc}
	    #1&#2 \\ #3&#4 \\#5&#6 \end{array} \right] $}}
\newcommand{\matrrrrcc}[8]{\mbox{$ \left[ \begin{array}{cc}
	    #1&#2\\#3&#4\\#5&#6\\#7&#8 \end{array} \right] $}}
\newcommand{\matrrrccc}[9]{\mbox{$ \left[ \begin{array}{ccc}                    
            #1&#2&#3\\#4&#5&#6 \\#7&#8&#9 \end{array} \right] $}}

\def\angb#1{\langle #1 \rangle }

\newcommand{\matrcc}[2]{\mbox{$ \left[ \begin{array}{cc} #1 & #2 \end{array}
           \right] $}}
\newcommand{\matrrc}[2]{\mbox{$ \left( \begin{array}{c}  #1 \\ #2 \end{array}
           \right) $}}
\newcommand{\matrrcc}[4]{\mbox{$ \left[ \begin{array}{cc} #1&#2 \\ #3&#4
           \end{array} \right] $}}
\newcommand{\IZ}{\Bbb{Z}}

\title[Grothendieck groups for operator algebras]{Grothendieck group
invariants for partly\\
self-adjoint operator algebras}
\author{S. C. Power}
\address{Dept.\ of Mathematics \& Statistics\\ Lancaster University\\
     Lancaster, U.K.\ LA1~4YF}
\email{s.power@@lancaster.ac.uk}
\thanks{22 January 1999}
\maketitle

\begin{abstract}
Partially ordered Grothendieck group invariants are introduced for
general operator algebras and these are used in the classification of
direct systems and direct limits of 
finite-dimensional complex incidence algebras
with common reduced digraph $H$  (systems of $H$-algebras).  
In particular the {\em dimension distribution
group} $G(A; C)$, defined for an operator algebra $A$ and a self-adjoint
subalgebra $C$, generalises both the $K_0$ group of a $\sigma -$unital
$C^*-$algebra $B$
and the 
spectrum (fundamental
relation)  $R(A)$ of a regular limit $A$ of triangular digraph algebras.
This invariant is more economical and computable
than the regular Grothendieck group $G^r_H(-)$ which nevertheless
forms the basis for a complete classification of regular systems
of $H$-algebras.

\end{abstract}

 \tableofcontents

\section{Introduction}

Two invariants have proven to be fundamental for the
classification of approximately finite-dimensional operator algebras. 
For self-adjoint algebras the scaled ordered $K_0$ group provides a
complete invariant for $C^*$-algebra isomorphism. See
Elliott \cite{ell}.
On the other hand for  triangular limit algebras determined
by regular embeddings
the so called spectrum, or fundamental relation, 
is a topological binary relation
which provides a complete invariant for isometric isomorphism. See
Power \cite{scp-book}.  In what follows we introduce various
Grothendieck group invariants, 
and in particular  
{\em dimension distribution groups}, which generalise these
and we
obtain classifications of systems and limits in terms of them.
\medskip

An operator algebra $A$ has associated with it a specific
inclusion $A \rightarrow {\cal L}(H)$ for some 
Hilbert space $H$ and so comes equipped with
a $C^*$-algebra inclusion $A \rightarrow C^*(A)$.  Accordingly
the natural homomorphisms between operator
algebras $A, A'$ are the star-extendible homomorphisms,
that is, those homomorphisms that 
are restrictions of $C^*$-algebra
homomorphisms $C^*(A) \rightarrow C^*(A')$.  The set
$Pisom (A)$ of partial isometries in $A$
{\em whose initial and final projections belong to A}
and the set $Pi(A) = Pisom
(A \otimes M_{\infty}$) associated with the stable algebra $A \otimes
M_{\infty}$  of $A$, are, as we shall see, a rich source of
invariants for star-extendible isomorphism.  
The dimension distribution groups given here 
are derived in this way
and in the case of approximately finite
algebras they provide information
on  the "rank distribution" of decompositions of
partial isometries
relative to the block structure of finite-dimensional superalgebras.
This is in parallel to the view that the $K_0$ class of a projection
determines its rank distribution.
\medskip

Recall that a (complex) {\em digraph algebra},  or finite-dimensional incidence
algebra, is a subalgebra $A$ of $M_n = M_n(\IC)$ which contains a 
maximal abelian self-adjoint subalgebra (a masa)  which, without loss of
generality, we may take to be the diagonal subalgebra.  In
Power \cite{scp-ko} (see also Heffernan \cite{hef}) we considered
the nonselfadjoint limit algebras
$ A  = \indlimit (A_k, \alpha_k) $
for which the embeddings $\alpha_k$ are star-extendible maps between
direct sums of $2 \times 2$ block upper triangular matrix algebras.  
This is the simplest family of partly self-adjoint algebras 
which are limits of digraph algebras and the classification may be effected
by the scaled $K_0$ group together with a binary relation on the scale
of $K_0$.
This binary relation, called the {\em algebraic order}, is determined
directly by
the set $Pisom{(A)}$.
\medskip

The case of limits of $r \times r$ block upper triangular matrices, for $r>2$,
is much more complicated. Star-extendible embeddings
between such matrix algebras are no longer automatically regular 
in the sense of being decomposable into a sum of multiplicity one embeddings
and even in
the case of regular embeddings, the multiplicity signature increases
rapidly with $r$.  Nevertheless, 
we obtain classifications of families of
regular direct systems $\ca = \{A_k, \alpha_k\}$
up to regular isomorphism.
For small values of $r$, or for restricted families of embeddings,
the induced group homomorphism 
$ G(\alpha_k) : G(A_k) \to G(A_{k+1})$ determines the inner unitary
equivalence class of $\alpha_k$ and we show how the limit group
$G(\ca)$ may be 
used in classification.
For more general families, such as the regular systems of $H$-algebras, more 
complicated Grothendieck group invariants, such as the 
scaled regular Grothendieck group
$G^r_H(\ca)$ are effective. This perspective is developed further
in \cite{scp-afa} where we use metrized variants of such ordered groups
in the classification of general limit algebras of irregular systems.
\medskip

For regular systems of digraph algebras we may view the
partial isometry dimension distribution group
quite naturally as the Elliott dimension group determined by the 
Bratteli diagram whose nodes represent 
the edges of the reduced digraphs of the building block algebras.
The improper edges, or loops, are
in correspondence with the vertices of the reduced digraphs and
so the sub-Bratteli diagram determined by these improper edges
gives rise to $K_0(\ca)$ and its inclusion in $G(\ca)$.
\medskip

The purely algebraic problem of classifying locally finite asociative algebras
has been emphasied recently by Zaleskii \cite{zal}.
The invariants given here seem to be
particularly useful in this consideration, at least for complex locally finite
associative algebras that are determined by regular embeddings.
In fact it seems likely that for most 
families of  regular digraph algebra systems,
and possibly for all such systems,  the 
star-extendible isomorphism class of the locally finite algebra
\[
A = \indlimit \ca
\]
determines the system $\ca$ up to regular star-extendible isomorphism.
This is known for a number of special cases; self-adjoint
systems, triangular systems, $T_2$-algebra systems,
systems of cycle algebras and various order preserving systems.
See Donsig \cite{don-cho} and Hopenwasser and Power \cite{hop-scp-2}
for example.

The paper is organised as follows. In Section 2 we define the 
dimension distribution group $G(\ca)$ for a regular system of digraph algebras
as well as the regular Grothendieck group invariants.
The dimension distribution group
$G(A; C)$ is defined for a nonselfadjoint operator algebra
$A$ containing a regular maximal abelian self-adjoint subalgebra
$C$. In the case of triangular AF algebras the group $G(A; C)$
coincides with the group of continuous integer valued functions
on the spectrum which vanishi at infinity, and both AF 
C*-algebras and TAF algebras are classified
by the same invariant.
In Section 3 the middle 
ground between these extremes is considered.
Regular direct systems 
of $T_2$-algebras and $T_3$-algebras are shown to be classified
in terms of the dimension distribution group and a natural
matrix unit scale.
Similar results hold for restricted classes of systems of $T_4$
algebras and for systems of cycle algebras.
In section 4 we indicate how regular Grothendieck groups are used in the
classification of certain regular systems of $H$-algebras and cycle algebras.

\section{Grothendieck group invariants}

\subsection{Regular Grothendieck groups}

First we introduce the general Grothendieck group invariants.
\medskip

Let $\ca = \{A_k, \alpha_k \}$ be a direct system of digraph algebras for       
which the algebra homomorphisms $\alpha_k: A_k \rightarrow A_{k+1}$
are star-extendible and regular.  Such systems will be the fundamental
objects of the sequel. Two systems $\ca, \ca'$ are said to
be {\em regularly isomorphic} if there exist star-extendible
regular maps $\phi_k, \psi_k$ forming a commuting diagram
\medskip
\[
            \begin{CD}  A_1 @>>>  A_{n_1} @>>>  A_{n_2} @>>> \cdots 
           \\
    & \lsymbse{\phi_1} & \symbup{\psi_1} & \symbse{\phi_2} & \symbup{\psi_2}
        & & \cdots  \\
A'_1 @>>>  A'_{m_1} @>>>  A'_{m_2} @>>> \cdots 
\end{CD}
\]
\noindent where the horizontal maps are those of the subsystems $\{A_{n_k} \},  
\{A'_{m_k} \}$ of $\ca, \ca'$ respectively.
\medskip
 
The commuting diagram above generally leads to many natural invariants
and the main aim in classifying
a given family is to find a minimal or economical set of
invariants which give complete invariants for regular isomorphism.
Thus, if $\cf$ is the family of {\it all} triangular systems then
the spectrum $R(\ca)$ provides a complete invariant.  On
the other hand, if $\ca$ is the subfamily of triangular
alternation algebras then  $R(\ca)$ may be replaced
by an equivalence class of a pair of
generalised integers. (See \cite{hop-scp}, \cite{poo}.)
Likewise in the partly self-adjoint
setting of rigid systems of 4-cycle algebras, which is indicated in
section 3.4, a certain jointly scaled group
$K_0(\ca)
\oplus H_1(\ca)$ provides complete invariants,
whereas for the subfamily of
even stationary matroid systems the abelian group
$H_1(\ca)$ together with the ordered $K_0$ group, suffices.
These facts may be found in \cite{don-scp-2}.
\medskip

Let $\cf_r, r = 1, 2, \dots, $ be the family of regular  
star-extendible systems with digraph algebras which are $T_r$-algebras,
where by a $T_r$-algebra we mean a block $r \times r$ upper triangular
matrix algebra.  Also,  for an antisymmetric transitive
digraph $H$, let $\cf_H$
denote the systems with digraph algebras $A_k$ whose reduced digraph is
$H$, for all $k$. We refer to such digraph algebras as
$H$-algebras.
\medskip

For a reduced digraph  $H$ define $V_H^r(\ca)$
to be the abelian semigroup
of equivalence classes $[\phi]$ of star-extendible regular
embeddings
\[
\phi : A(H) \to \ca \otimes M_\infty.
\]

\noindent That is we require $ \phi $ to be a regular 
star-extendible embedding from
 $ A(H)$ to $ A_k \otimes M_n$  
for some $n$ and some algebra $A_k$ of the direct system $\ca$,
and such embeddings $\phi, \psi$ are equivalent if there exists
a unitary $u$ in
$A_j \otimes M_N$
for some $j, N$ such that
$ u\phi(a)u^* = \psi(a)$ and $ \phi(a) = u^*\psi(a)u$ for all $a$ in $A(H)$.
The additive semigroup operation is given by
$[\phi] + [\psi] = [\phi \oplus \psi]$.
Define $G^r_H(\ca)$ to be the Grothendieck group
of $V_H^r(\ca)$. It can be shown that $V_H^r(\ca)$ has cancellation
and so embeds injectively in $G^r_H(\ca)$.
We refer to $G^r_H(\ca)$ as the {\it regular Grothendieck group
invariant} of $\ca$, for the digraph $H$, the terminology
reflecting the fact that there are no restrictions
on the embeddings beyond regularity and star extendibility. 

Note that if $H$ is the trivial digraph with one vertex and edge
then $G^r_H(\ca)$ is naturally isomorphic 
to $K_0(\ca)$ through the correspondence
of $[\phi]$ with the $K_0$ class of the projection $\phi(1)$.

On the other hand if $A(H) $ is the upper triangular
matrix algebra $ T_4$ then, as we see later, there are
35 classes of multiplicity one embeddings from $T_4$ to $T_4 \otimes M_\infty$.
It follows that the regular Grothendieck group of $T_4 \otimes M_k$
for the digraph for $T_4$
is isomorphic to ${\Bbb Z}^{35}$.
Thus if $\ca$ is a regular system of $T_4$-algebras
then $G^r_{T_4}(\ca)$ is a dimension
group of the form  $\indlimit ({\Bbb Z}^{35}, \hat{\alpha_k})$ 
where each $\hat{\alpha_k}$
is realised by a nonnegative integral
matrix.
\medskip

The regular Grothendieck group  $G^r_H(\ca)$  admits a natural ordering
which is induced, as is the case for $K_0$, by a
scale $\Sigma_H^r(\ca) \subseteq G^r_H(\ca)$.
This scale is defined to be the image in $G^r_H(\ca)$ of the
classes of those embeddings $\phi$ which map into
$\ca$ rather than the stable system $\ca \otimes M_\infty$. The
partially ordered scaled abelian group $(G^r_H(\ca), \Sigma_H^r(\ca)) $ 
is clearly an invariant for regular isomorphism
and in Section 4 we indicate how the pair
$(G^r_H(\ca), K_0(\ca))$ 
may be used in the classification of
the regular systems for  $\cf_H$.

First we turn to the consideration of a more computable invariant.

\subsection{\bf The dimension distribution group $G(\ca)$}

For an $AF \; C^*$-algebra $B = \indlimit \ B_k$
the  $K_0$ class of a projection $p$ in $B_l$ 
determines how the rank of $p$
is distributed in the summands of $B_k$ for
$k > l$.  The idea behind the formulation of the group $G(A; C)$ is to
capture in a similar way how the rank of a partial isometry in a
non-self-adjoint limit algebra $A = \indlimit \ A_k$ 
is distributed
in the summands and in the block subspaces of $A_k$ for large $k$. 
\medskip

Let $\alpha: T_3 \otimes M_m \rightarrow T_3 \otimes M_n$,  with $n$
larger than $m$, be a star-extendible algebra 
homomorphism.
Suppose,  furthermore,  that the image of each standard matrix unit
$e_{ij} \otimes f_{kl}, \; \; 1 \leq i \leq j \leq 3, \; \; 1 \leq k,
l \leq m$, is a
sum of standard matrix units.  Then to each image
\[
\alpha (e_{ij} \otimes f_{kl}) =
\matrrrccc{v_{11}}{v_{12}}{v_{13}}{0}
{v_{22}}{v_{23}}{0}{0}{v_{33}}
\]
we may associate the {\em rank distribution matrix}  in
$T_3(\ZZ)$ defined by 
\[
rk(\alpha (e_{ij} \otimes f_{kl})) = 
\matrrrccc{\mbox{rank}(v_{11})}{\mbox{rank}(v_{12})}
{\mbox{rank}(v_{13})}{0}{\mbox{rank}(v_{22})}
{\mbox{rank}(v_{23})}{0}{0}{\mbox{rank}(v_{33})}
\]
and this matrix is independent of $k,l$.  The correspondence
$$rk(e_{ij} \otimes f_{11}) \rightarrow rk(\alpha(e_{ij} \otimes
f_{11}))$$
extends to an abelian group homomorphism $\hat{\alpha} : T_3({\ZZ})
\rightarrow
T_3({\ZZ})$. Moreover, under the natural identification $K_0(T_3 \otimes M_m) 
$ is equal to $
D_3 ({\ZZ})$, the diagonal subgroup, and the restriction
$\hat{\alpha}|D_3({\ZZ})$ coincides with $K_0 \alpha$.  Plainly, 
$(\beta \circ \alpha)^{\hat{}}
=
\hat{\beta } \circ  \hat{\alpha }.$
\medskip

More generally we may define  $\hat{\alpha}$ in a similar way when 
$\alpha : A_1 \to A_2$ is a {\em regular} star-extendible 
homomorphism between digraph algebras,
that is,  one which is a direct sum of multiplicity one
embeddings. In this case
matrix units can be chosen for $A_1 $ and $A_2$
so that matrix units map to sums of matrix units.
\medskip

Thus, given a direct system
$\ca = \{A_k, \alpha_k \},$
we can 
introduce the abelian
group
\[ 
G(\ca) = \lim_{\rightarrow} (G(A_k), \hat{\alpha}_k).
\]
together with, as we see below, a natural ordering and a scale
$\Sigma (\ca)$ which augments the $K_0$ ordering and scale.  
We call $G(\ca)$ the (scaled ordered) {\em dimension
distribution group}, or the partial isometry dimension distribution group,
of the system $\ca$.
In fact there are scaled ordered group 
surjectons 
\[
\pi_f : G(\ca) \to K_0(\ca), \pi_i : G(\ca) \to K_0(\ca)
\]
and the quadruple 
\[
((G(\ca), \Sigma(\ca)), (K_0(\ca), \Sigma_0(\ca), \pi_f , \pi_i)
\]
may  be viewed as a generalisation of a dual form of
the spectrum of a triangular 
limit algebra.
 
\subsection{\bf The dimension distribution group $G(A; C)$}

Fix a  general operator algebra $A$ and 
a self-adjoint subalgebra $C$.  
We now give a general formulation of a 
dimension distribution group invariant $G(A; C)$.
However our principal concern is for a pair
$(A, C)$ where
\[
 A=\indlimit \  \ca = \indlimit \ (A_k, \alpha_k),\ 
\ \ \ \ C = \indlimit \  \cc = \indlimit \  (C_k, \alpha_k)
\]
are associated with a regular star-extendible digraph 
algebra system $\ca$
and  subsystem $\cc$ of masas
$C_1, C_2, \dots$ 
such that, for all $k$, $\alpha_k$ maps the partial isometry normaliser
of $C_k$ into the partial isometry normaliser
of $C_{k+1}$. 
Recall that an element $v$ in a star algebra $B$ normalises
a subalgebra $E$ of $B$ if $vEv^* \subseteq E$ and $v^*Ev \subseteq E$.
\medskip

Write $M_\infty(A)$ for the usual stable algebra $M_\infty({\IC}) \otimes A$
of $A$, write
$D_{\infty}(A)$ for its diagonal subalgebra and 
define the following sets of projections and partial isometries.
\medskip

$$Pr(A) = \{p \ \in \  M_{\infty}(A): p =p^2 = p^* \},$$
$$Pi(A) = \{v \ \in \  M_{\infty}(A): v^*v, vv^* \ \in \  Pr (A)
\},$$
$$\tilde{P}r(A; C) = \{p \ \in \  Pr (A): p \mbox{ normalises } 
D_{\infty}(C) \},$$
$${P}r(A; C) = \{p \ \in \  Pr(A): p \mbox{ normalises } M_{\infty}(C) \},$$
$$\tilde{P}i(A: C) = \{v \ \in \  Pi(A) : v  \mbox{ normalises } D_{\infty}(C) \},$$
$$Pi(A; C) = \{v \ \in \  Pi(A): v \mbox{ normalises } M_{\infty}(C) \}.$$
\medskip

Also write $Proj(A)$, $Pisom(A)$, $Proj (A; C)$, $Pisom (A; C)$ for the
appropriate sets in $A$, rather than the stable algebra of $A$.  In 
particular $Pisom (A)$ is the set of partial isometries in $A$
{\em with initial and final projections in $A$}, and $Pisom (A;C)$ 
is the subset of elements which
normalise the masa $C$.
\medskip

For $v, w$ in $Pi(A)$ write $v \simeq w$ if $zvz^* = w$ for some unitary
$z$ in $M_{\infty}(A)$. More precisely $zvz^* = w$ with $z$  a unitary 
in some sufficiently large matrix algebra containing $v$ and $w$.
Write $v \sim w$ if $zvy = w$ for some
unitaries $z, y$ in $M_{\infty}(A)$.  Similarly for the four normalising
sets above, define unitary equivalence and equivalence with respect to
the appropriate normalising unitaries.  The resulting 
equivalence classes provide the twelve
monoids,  
\[
Pr(A)/\sim, \ \ \  Pr(A)/\approx, \ \ \dots \ \ ,\ Pi(A;C)/\approx \  
\]
and twelve
associated Grothendieck groups.  
\medskip

In the case of finite-dimensional digraph algebras one may also consider
the following subset of $Pisom (A)$ consisting of (let us say) {\em
regular partial isometries}:
\[
Pisom_{reg}(A) = \{v \in Pisom(A) : v = v_1 + \dots +v_k, 
v_i
\ \in \  Pisom (A), v_i \mbox{ rank 1 } \}.
\]
Then it is elementary to obtain the identifications
\[
Pisom_{reg} (A)/\sim \ \  = \ \ Pisom (A; C)/\sim,\ \ \ 
Pi_{reg}(A)/\sim \ \ =
\ \ Pi(A; C)/\sim,
\]
in this setting.
\medskip

Returning to the general pair $A, C$, note that
for projections $p,q$ if $zpy=q$, with $z,y$ unitary, then $y^*py = q$. 
Thus, 
\[
Pr(A)/\sim \ \ = \ \ Pr(A)/\approx,
\]
 with similar equalities for
$Pr(A; C)$
and $\tilde{P}r(A; C)$.  Since, in general, 
\  $\sim$ equivalence of partial isometries extends
unitary equivalence of projections it is natural to
restrict attention to the
following Grothendieck groups:
\medskip

$$K_0(A) = Gr (Pr (A)/\approx),$$
$$K_0(A; C) := Gr(Pr (A; C)/\approx),$$
$$\tilde{K}_0(A;C) := Gr (\tilde{P}r(A;C)/\approx),$$
$$G_t(A) := Gr(Pi(A)/\sim),$$
$$G(A;C) := Gr(Pi(A;C)/\sim),$$
$$\tilde{G}(A;C) := Gr(\tilde{P}i(A;C)/\sim).$$
\medskip

These groups are appropriate in the context of $\sigma -$unital
local operator algebras, where by {\em local} we mean, as in 
Blackadar \cite{bla},
that $A$ (and $C)$ are closed under holomorphic functional calculus. 
In particular, an unclosed union of closed subalgebras is a local operator
algebra, as are algebraic direct limits of norm-closed operator algebras.
 In this context $Gr (Pr(A)/\approx)$ is the usual $K_0$ group of $A$. 
Also if $B$ is a $\sigma$-unital $C^*$-algebra then
$$G(B;B) = K_0(B).$$
\medskip

The group $G_t(A)$, which we refer to as 
the "total partial isometry" dimension distribution group,
does not yield a convenient invariant for digraph algebras and their
limit algebras.
At least, $G_t(A)$ seems inappropriate for the class of 
regular limit algebras.
This is indicated below in Example 2.
 Accordingly we make the following definition of $G(A)$
in this case, and with this definition we have the natural identification
$G(\ca) = \indlimit \  G(A_k)$ for a star-extendible regular 
system of digraph algebras.
\medskip

\begin{defn}
  The dimension distribution group $G(A)$ of a digraph 
algebra $A$ is defined to be the group $G(A; C)$, where $C$ is a masa in $A$.
\end{defn}
\medskip

\noindent{\bf Example 1. } $G_t(T_2({\IC}))$ 
is the abelian group $ T_2({\ZZ}).$  In fact if $v
\ \in \  Pisom(T_2 \otimes M_n)$ then
\[
v =  \left[ \begin{array}{cc}
            v_1& v_2 \\ 0 &v_3  \end{array} \right] 
\]
with respect to the block upper triangular decomposition. From the
definition of $Pisom(~)$ it follows that 
$v^*v$ and $vv^*$ are block diagonal.  Thus 
$v_1$ and $v_3$, and hence $v_2$ are partial isometries.
 The correspondence $[v] \rightarrow rk(v)$ leads to the stated
identification.  
\medskip

\noindent{\bf Example 2.} If $n \geq 3$ then $G_t(T_n({\IC}))$ contains ${\ZZ}^c$, where $c$ is the
cardinality of the continuum.  Indeed, there are
uncountably many $\sim$ equivalence classes in $T_3({\IC})$ of partial
isometries of the form
\[
\matrrrccc{0}{v_1}{v_2}{0}{v_4}{v_3}{0}{0}{0}
\]
where the $2 \times 2$ submatrix is
a unitary matrix in $M_2$. These classes lead to an injection
${\ZZ}^c_+ \rightarrow Pi(T_3({\IC}))/\sim$.  The proliferation of
such equivalence classes is one 
reason for the consideration of normalising
partial isometries.
\medskip

\noindent{\bf Example 3.}  For $n \geq 0 $,
$ G(T_n; D_n) = \tilde{G}(T_n; D_n) = T_n({\ZZ})$.
More generally, if $X$ is a totally disconnected locally compact
Hausdorff space, then 
$$G(T_n \otimes C(X); D_n \otimes C(X)) = \tilde{G}(T_n \otimes C(X); D_n
\otimes C(X))$$
$$= T_n({\ZZ}) \otimes K_0(C(X))$$
$$ = {\ZZ}^{n(n+ 1)/2} \otimes C(X, {\ZZ})$$  .
\medskip

One can establish this identification and that of the next example,
by a simple direct
argument, or by appealing to the more general discussion of the next section.
\medskip

\noindent{\bf Example 4.}
Let $A = \indlimit (T_{2^k}, \rho_k)$ be the $2^\infty$ refinement algebra
with spectrum $R(A)$. (See \cite{scp-book} and below.) 
Let $B$ be a UHF C*-algebra
with $K_0$ group ${\Bbb Q}_B $ realised as an ordered unital subgroup 
of $  {\IQ}$. Also let $C$ be a masa of the form
$(A\cap A^*) \otimes D$ with $D$ a regular canonical masa of $B$. Then
\[
G(A \otimes B, C) = C(R(A), {\Bbb Q}_B),
\]
where the right hand side is the group of continuous 
${\Bbb Q}_B$-valued functions
on the spectrum which vanish at infinity.
In fact the masa $C$ is intrinsic to the algebra $A \otimes B$ as the unique 
(up to star-extendible automorphism) standard AF diagonal. See Theorem 4.1 of
\cite{scp-hom2}.
Thus we may write $G(A \otimes B)$ for $G(A \otimes B; C)$.

\medskip

\subsection{\bf $G(A; C)$ and the spectrum $R(A)$}

Let $A = \indlimit \  (A_k, \alpha_k)$ be a limit of digraph
algebras with respect to regular star-extendible embeddings.  Then for
$k=1, 2, \dots$ there are matrix unit systems $\{e^k_{ij} \}$ for the
finite-dimensional $C^*$-algebras $B_k = C^*(A_k)$ such that $A_k$ is
spanned by some of these matrix units, including the diagonal matrix
units, and such that $\alpha_k$ maps matrix units in $\{e^k_{ij}\}$ to
sums of matrix units in $\{e^{k+1}_{i,j} \}$.  Thus if $C_k$ is the masa
of $A_k$ (and $B_k$) spanned by $\{e^k_{ii}\}$ then we have the triple
$C \subseteq  A \subseteq B$, where 
\[
B = C^*(A) = \lim_{\rightarrow} (B_k, \tilde{\alpha}_k),
\]
with $\tilde{\alpha}_k$ the unique extension of $\alpha_k$ and where $C$ is the
{\it regular canonical masa}
 \[C = \lim_{\rightarrow}(C_k, \alpha_k).\] 
\medskip

We now recall the definition of the fundamental relation or { topological
binary relation}, or {\em spectrum}, $R(A;C)$ for the pair $(A,C)$.  See
\cite{scp-class}, \cite{scp-tensor} and Chapter 7 of \cite{scp-book}. 

\medskip
  Each matrix unit $e= e^k_{ij}$ induces a partial homeomorphism
$\alpha$ of the Gelfand space $X$ of $C$  such that 
$$\alpha(y)(ece^*) = y(c)$$
for all characters $y$ with $y(e^*e)=1$.  The set of such $y$ is the
domain $d(\alpha)$ of $\alpha$ and   the graph of $\alpha$ is the subset $E$
of $X \times X$ given by $E = \{(x,y): y \ \in \  d(\alpha),
x=\alpha(y)\}.$
\medskip

The binary relation $R(A; C)$ is then defined as
\[R(A; C)= 
\bigcup_{i,j,k} \{E^k_{ij}: E^{k}_{ij} \mbox{ is the graph for } e^k_{ij}
\ \in \  A_k \}\]
and the topology on $R(A; C)$ is the smallest topology for which each
$E^k_{ij}$ is both closed and open.  
In particular $R(A; C)$ is a locally compact topological space.
It can be shown that $R(A; C)$ is
well-defined and independent of the particular choice of direct system
or matrix units for the pair $(A,C)$. See \cite{scp-book}.
\medskip

When $A$ is triangular, that is, when $A \cap A^* = C$, then 
$A$ is referred to as a TAF algebra and we 
write $R(A)$ for $R(A; C)$. This topological
binary relation is a complete invariant for isometric isomorphism
and has been a useful tool for classification and analysis.
This 
can be seen, for example, in various
articles which have appeared since \cite{scp-book}, such as
\cite{hop-lau}, \cite{hop-pet}, \cite{poo-wag},
 \cite{scp-lex2},
 \cite{don-hop}, \cite{don-hud},
\cite{pet-poo-lex}.
\medskip
 
Let us now consider the dimension distribution group $G(A;C)$ with $A$
as above and assume that $A$ is 
a TAF algebra.  
\medskip

\begin{prop}
If $A$ is a TAF algebra then $G(A; C) = C(R(A), {\Bbb Z})$.
\end{prop}

\begin{pf}
Note first that the abelian semigroup $C(R(A),
{\ZZ}_+)$ of positive integer-valued continuous functions on the spectrum
$R(A)$  (vanishing at infinity will be understood) is 
generated by the family of characteristic functions $\chi_E$ for
closed-open sets $E$ associated, as above, with the matrix units $e$ of
the fixed matrix unit system $\{e^k_{ij} \}$.  It is a fundamental fact,
which follows quickly from Lemma 5.5 of \cite{scp-book}, that if $v
\ \in \  Pisom (A;C)$ then $v$ is unitarily equivalent by a unitary
in $C$ to a partial isometry which is a finite sum of matrix units. 
Furthermore, note that if $v, w \ \in \ Pisom (A; C)$ then, since
$A$ is triangular, $v \sim w$ if and only if $v=cw$ for some element
$c$ in $C$.  Thus to the class $[v]_{\sim}$ in $Pisom (A; C)/\sim$ there is a
unique element $z$ in $Pisom (A; C)$ which is a sum $z= e_1 + \dots +
e_r$ of the given matrix units, such that $[v]_{\sim} = [z]_{\sim} $.  
Although the
decomposition $z = e_1 + \dots + e_r$ is not unique the associated
correspondence
$$\kappa : [v] \rightarrow \chi_{E1} + \dots + \chi_{Er}$$
gives a well-defined semigroup injection
\[
Pisom(A; C)/\sim \ \ \ \rightarrow \ \ \ C(R(A), {\ZZ}_+).
\]
\medskip
 
Consider now $v \ \in \  \tilde{P}i(A; C)$ so that $v \ \in \  
M_n(A) \subseteq
M_{\infty}(A)$, for some $n$, and $v$ normalises $D_n(C)$.  By the
reasoning above, with $D_n(C)$ playing the role of $C,$  $v$ is
unitarily equivalent to a partial isometry $e \ \in \  M_n(A)$ of
the form $e=(g_{ij})_{i,j=1}^n $ where each $g_{ij}$ is a sum of matrix
units for the given matrix unit system for $A$.  Since $e$ 
normalises $D_n(C)$ the partial isometries
$\{g_{ij}\}_{i,j=1}$ are orthogonal, that is, they have orthogonal range
projections and orthogonal initial projections.  It follows
that in $M_{n^2}(A)$ we have $e \sim w$ where $w$ is the block diagonal
partial isometry
$$w = \sum^{n}_{i=1} \sum^{n}_{j=1} \oplus g_{ij}.$$
If we define 
$$\kappa: Pi(A;C)/\sim \ \ \rightarrow \ \ C(R(A), {\ZZ}_+)$$
by
$${\kappa}([w]_\sim) = \sum^{n}_{i=1} \sum^{n}_{j=1} \kappa([g_{ij}]_\sim).$$
then $\kappa$ is a well defined semigroup surjection. We claim that $\kappa$
is injective. We have, by orthogonality
again, $e \sim f$ where $f$ is a block diagonal sum \[f= \sum^N_{j=1}
\oplus
f_j\] in $D_N(A)$ with each $f_j$ a matrix unit in $A_M$, for some
sufficiently large $M,N$.  If \[f'= \sum^N_{j=1} \oplus f_j'\] is
another
such sum, for the same $M, N$, but with, perhaps, some summands equal to zero,
then $f \sim f'$ if and only if each $f_j^\prime$ in the
sum for $f^\prime$ appears in the sum for $f$, with the same multiplicity,
and vice versa.
However, this is precisely the condition that
$$\sum^N_{j=1} \chi_{F_j} = \sum^N_{j=1} \chi_{F_j^\prime},$$
and so it follows that $\kappa$ is injective.
\medskip

Thus 
$$\tilde{G}(A; C)= Gr(\tilde{P}i(A;C)/\sim)= 
Gr(C(R(A), {\ZZ}_+)) = C(R(A), {\ZZ}).$$
\medskip

We have shown, in Lemma 2.4 of \cite{scp-hom2},
that a partial isometry $w$ in $M_n \otimes B$ which normalises 
$M_n \otimes C$ has the form  $dw_1$ where $d$ is a partial isometry 
in $M_n \otimes C$ and $w_1$ is a partial isometry in $B_k$,
for some $k$, which normalises $C_k$. From this one can deduce that
there is a natural identification of $G(A; C)$ with $\tilde{G}(A; C)$.
Thus
$$G(A; C) = C(R(A), {\ZZ}).$$
\end{pf}
\medskip

At the other extreme we have
\medskip

\begin{prop}
Let $C$ be a standard AF diagonal subalgebra in the 
AF C*-algebra $B$ (as above).
Then 
$$G(B; C)= K_0(B)$$
\end{prop}
\medskip

\begin{pf}
The identification follows from the fact that
normalising partial isometries in a
finite-dimensional $C^*$-algebra are $\sim$ equivalent, by a normalising
unitary, if and only if they have the same rank distribution in the
summands of the $C^*$-algebra.
\end{pf}
\medskip

\subsection{\bf Orders and Scales}

\begin{defn}
The {\em positive cone} $G(A;C)_+$ of $G(A;C)$ is defined to be the
semigroup generated by the (images of) classes $[v]$ where $v
\ \in \  Pisom (A;C)$ and the {\em scale} of $G(A;C)$ denoted $\Sigma (A;C)$
is defined to be
the set of these (images of) classes.
\end{defn}
\medskip

In particular  $G(T_n; D_n)_+ = T_n({\ZZ}_+)$ and $\Sigma (T_n; D_n)$ 
is the set of
zero-one elements with at most one nonzero entry in each row and column.
This follows from the next
proposition.
\medskip

Let $A$ be an algebra with 
$D_n \subseteq A \subseteq M_n$,
that is, let  $A$ be a digraph algebra.
Then there are natural scaled
group homomorphisms
$$\pi_f: G(A) \rightarrow K_0A, $$
$$\pi_i: G(A) \rightarrow K_0A$$
which extend the correspondences
$$[v] \rightarrow [vv^*], \ \ [v] \rightarrow [v^*v]$$
for $v \ \in \  Pi(A; D_n)$.  We may assume without loss of
generality that 
\[
A = span\{e_{ij}:e_{ij} \ \in \  E(G) \}
\]
where
$\{e_{ij} \}$ is the standard matrix unit system, and $E(G)$ is the set
of edges of a directed graph with vertices $\{1, \dots, n\}$.  
Let $G_r$ be
the {\em reduced graph} of $G$ 
obtained from $G$ by
identifying vertices which are equivalent.
(Thus, vertex $i$ is equivalent to vertex $j$ if the
matrix units $e_{ii}, e_{jj}$  lie in the same matrix algebra
summand of $A \cap
A^*$.)
Then 
$G(A)={\ZZ}^{|E(G_r)|}$.  
 The digraph algebra $A(G_r)$ is a triangular matrix
algebra, with matrix units $f_{kl}$ say, and there is a natural
identification of $G(A)$ with the abelian (additive) subgroup
generated by the matrix units $\{f_{kl}: (k,l) \ \in \  E(G_r) \}.$
A straightforward
 induction proof yields the following connection between scales.

\medskip

\begin{prop}
  {\it Let $A$ be a digraph algebra and let $g \ \in \  G(A)_+.$   
Then $g
\ \in \  \Sigma (A;D_n)_+$ if and only if $\pi_i(g)$ and $\pi_f(g)$
belong to $\Sigma_0 A$, the scale of $K_0A$. }
\end{prop}
\medskip

\begin{defn}
  (i) Let $A = {\indlimit}
(A_k, \alpha_k)$  
be a regular star-extendible direct system of digraph algebras with
regular canonical masa $C=  \indlimit \ (C_k, \alpha_k)$. 
Then the
dimension distribution group invariant, or the {\em dimension distribution
group quadruple}, for the pair $A, C$ is the quadruple 
\[\cg(A; C)
= (G(A; C), K_0(A), \pi_f, \pi_i).
\]
\medskip

(ii) Two such quadruples, for $A, C$ and $A', C'$, respectively,
are said to be isomorphic if
there is a $\pi$-respecting scaled group 
isomorphism $\theta: G(A; C) \rightarrow
G(A';C')$, that is, if there are two commuting diagrams
$$
\begin{CD}
G(A; C) @>{\theta}>> G(A';C')
\\
@V{\pi}VV @V{{\pi '}}VV  \\
K_0A @>{\mathrm{\theta | K_0(A)}}>> K_0A'
\end{CD}
$$
where the horizontal maps are scaled group isomorphisms
and where $\pi = \pi_i$ or $\pi_f$.
\end{defn}
\medskip

Of course the distribution group quadruple is indeed an invariant for
the pair $(A,C)$, with respect to masa preserving star-extendible
isomorphism.  In the self-adjoint case it degenerates to the scaled
$K_0$ group, and so, by Elliott's theorem $\cg(A; C)$ is a well defined
complete invariant for 
AF C*-algebras. 
At the other extreme,
in the triangular case,  where there is only one
distinguished masa, 
the distribution group quadruple
is an invariant for star-extendible isomorphism. In fact
it is a complete invariant in this case also.

\begin{thm}
  { Let $A, A'$ be limit algebras of
 regular star-extendible systems of digraph algebras
and suppose that $A, A'$ are either self-adjoint or triangular.
Then $A, A'$ are
star extendibly isomorphic if and only if their dimension distribution
group quadruples  are isomorphic. }
\end{thm}
\medskip

\begin{pf}
 We have already commented on the self-adjoint case
so assume that  $A,A'$ are triangular and their 
quadruple invariants are isomorphic.
By Proposition 2.3 there are two
commuting diagrams
$$
\begin{CD}                                                                      
C(R(A), \Bbb{Z}) @>{\theta}>> C(R(A'), \Bbb{Z}) 
\\
@V{\pi}VV @V{{\pi '}}VV  \\
C(X, \Bbb{Z}) @>{\mathrm{\theta | K_0(A)}}>> C(X', \Bbb{Z})
\end{CD}
$$
with $\pi = \pi_f$ for one diagram and 
$\pi = \pi_i$ for the other.
Thus 
the isomorphisms
$\theta$ and $\theta|K_0 A$ are induced by homeomorphisms $ \nu : R(A)
\rightarrow R(A')$ and $\gamma : X \rightarrow X'$, respectively.  
In
view of the commuting diagrams,
it follows that $\nu = \gamma \times 
\gamma$.  Thus $\nu $ is a topological binary relation isomorphism.
Since the spectrum is a complete invariant (by \cite{scp-ko})
it follows that $A$ and $A'$ are star-extendibly
isomorphic.  
\end{pf}

Note that the scale  of $G(A; C)$ plays a crucial
role in the self-adjoint context and a rather 
redundant role in the triangular
context. 
\medskip

Recall that the {\it algebraic order} $S(A; C)$ for the pair
$(A,C)$ is the transitive binary relation in $\Sigma_0(A) \times \Sigma_0(A)$
which is the range of the correspondence
\[
[v] \rightarrow ([vv^*], [v^*v])
\]
for $v \ \in \  Pisom(A; C)$. 
In the next section we shall see that it is also fruitful
to consider more elaborate "matrix unit scales" in higher products
$\Sigma_0(A) \times \dots \times \Sigma_0(A)$.
\medskip

\section{Classifications with $G(\ca)$}

We now turn to the 
classification of specific regular
systems.
\medskip

Given an
isomorphism of invariants
there are two fundamental issues to address when  constructing a
regular isomorphism between systems $\ca$ and $\ca'$.
These are, firstly, the existence
of digraph algebra homomorphisms, with a given correspondence of
invariants, and, secondly, the uniqueness up to inner conjugacy of
these same homomorphisms.  Resolving the existence issue enables
the construction of $\phi_1$ and an initial $\psi_1$, as
in the diagram in Section 1, whilst
uniqueness enables the correction of $\psi_1$ to obtain a commuting
triangle.  
It is usually the existence question which is a more subtle issue.
\medskip

The dimension distribution group $G(\ca)$ can equally well be defined
for systems of digraph spaces and regular {\em linear} maps.
This suggests that the lifting of a restriction of a dimension
distribution group isomorphism to digraph algebra
{\em homomorphism} is conditional on its preservation of an
invariant which reflects the mutiplication in some way.  
\medskip

\subsection{Matrix unit scales.}

Consider  a $\pi_i, \pi_f$ respecting ordered group homomorphism 
$$\theta: G(T_3) \rightarrow G(B)$$
where $B$ is a digraph algebra.
Plainly $\theta$ is liftable to a 
regular star-extendible homomorphism if and only if
the triple 
$$\{ \theta(e_{12}), \theta(e_{23}), \theta(e_{13}) \}$$
belongs to a {"matrix unit scale"} 
in $( \Sigma (B) \times \Sigma (B) \times \Sigma (B)$ given by
$$\{([\phi(e_{12})], [\phi(e_{23})], [\phi(e_{13})]): \phi: T_3
\rightarrow B \}$$
where $\phi$ runs over star-extendible regular homomorphisms. 
In fact the property of $\pi_i, \pi_f$ respecting could be incorporated into
matrix unit scale preservation by enlarging the scale
to the set
\[
\{\Pi_{1\le i \le j \le 3} [\phi(e_{i,j}]:  \phi : T_3
\rightarrow B \} \subseteq (\Sigma(B))^{(6)}.
\]
We  refer to this ordered 
set as the {\it matrix unit scale} or, more precisely, 
as the $T_3$-matrix unit scale.
Thus, as we see in the  existence lemma below,
one way to ensure existence, that is, to ensure the
liftability of a restriction $\Phi|G(A_1)$ of a given scaled group
isomorphism $\Phi: G(\ca) \rightarrow G(\ca')$, is to demand that $\Phi$
preserve the appropriate matrix unit scale. 
\medskip

\begin{defn}
(i) Let $H$ be an antisymmetric digraph and let $B$ be a digraph algebra.
Then the 
{\em $H$-matrix unit scale} $S_{H} (B)$ for $B$ is the subset of
\[ 
\Pi_{(i,j) \varepsilon E(H)} \; \; \Sigma (B)
\]
given by the set of elements of the form 
\[
\Pi_{(i,j) \varepsilon E(H)} [\phi(e_{ij})]
\]
where  $(i,j)$
ranges through the edges of $H$ and where $\phi$ ranges over regular
star-extendible homomorphisms from $A(H)$ into $B$.

~~~(ii) If $\ca$ is a regular star-extendible system $\{A_k, \phi_k\}$
of  digraph algebras
then the {\it $H$-matrix unit scale}
 $S_{H}(\ca)$ 
is defined to be the union of the images of
the scales $S_{H}(A_k)$ in the ordered product
$\Sigma (\ca))^{\nu}$, where $\nu$
is the number
of edges of $H$.
\end{defn}
\medskip

Note that the inclusion $S_H(A_k) \to S_H(A_{k+1})$
is simply the direct product of the coordinate inclusions.
\medskip

\begin{lma}  
 \ \ \mbox{(Existence lemma.)} 
Let $A_1, A_2$ be $H$-algebras and let $\theta : G(A_1) \to G(A_2) $
be a group homomorphism. 
Then the following conditions are equivalent.
\medskip

(i) There is a regular star-extendible homomorphism
$\phi : A_1 \to A_2$ with $G(\phi) = \theta.$
\medskip

(ii)
$\theta(S_{H}(A_1)) \subseteq S_{H}(A_2)$.
\end{lma}   

\medskip

\begin{pf}
 Assume that $H_1$ and $H_2$ are the digraphs of 
$A_1, A_2$ so that the reduced digraphs of $H_1$ and $H_2$ 
are equal to $H$.
Let $i: H \to H_1$ be a digraph injection with
associated (multiplicity one) injection
$\alpha : A(H) \to A(H_1)$ for which 
the induced group homomorphism $G(\alpha)$
is the identity. With this injection we  identify $A(H)$ 
as a subalgebra of $A_1$. Plainly each matrix
unit $e$ in $A_1$ admits a unique 
representation $e = e_1fe_2$ with $f$ a matrix 
unit for $A(H)$ and $e_1, e_2$ matrix units
for $A_1 \cap A_1^*$. Let
\[
\beta_0 : A_1 \cap A_1^* \to A_2
\]
\[
\beta_1 : A(H) \to A_2
\]
be homomorphisms which map matrix units to 
sums of matrix units, and for which $K_0\beta_0 = \theta |\Sigma_0A_1$,
and $G(\beta_1) = \theta$.
The existence of the C*-algebra homomorphism
$\beta_0$ follows from the inclusion 
$\theta(\Sigma_0 A_1) \subseteq \Sigma_0A_2$, which is implied by
condition  (ii).
The existence of the star-extendible  homomorphism $\beta_1$
also follows from the condition (ii).
In order to define $\beta_1$
we only need limited information from the second condition,
namely that
\[
\Pi_{(i,j) \varepsilon E(H)} \theta(g_{ij}) = 
\theta(\Pi_{(i,j) \varepsilon E(H)} (g_{ij}) \in S_{H}(A_2),
\]
where here we write
$g_{ij}$ for the natural basis elements of $G(A(H))$, indexed by edges of $H$.
Having specified $\beta_1$ we may choose  $\beta_0$ 
in such a way that these maps agree on
$A(H) \cap A(H)^*$.
By the unique  factorisation of matrix units mentioned 
above there is a well-defined
linear map $\phi : A_1 \to A_2$ which extends 
both $\beta_1$ and $\beta_0$ and which satisfies
$\phi(e_1fe_2) = \phi(e_1)\phi(f)\phi(e_2)$. 
It follows that $\phi$ is an algebra homomorphism, and,
since $\beta_1$ is star-extendible, it follows that $\phi$ is star-extendible.

\end{pf}  
\medskip

One can also obtain an analogue of the lemma for a restricted family 
of embeddings.

\begin{defn}
Let $\Omega$ be a family of regular star-extendible embeddings between
$H-algebras$ which is closed under compositions and direct sums, and let
$\cf (\Omega)$ be the corresponding family of systems $\ca$.
 
~~~(i) If $A$ is a digraph algebra then the $\Omega$-scale, $S_\Omega (A)$,
is the subset of $S_{H}(A)$ corresponding to
homomorphisms
$\phi : A(H) \to A$ which belong to $\Omega$.
 
~~~(ii) If $\ca \in \cf (\Omega)$ then the $\Omega$-scale, $S_\Omega (\ca)$,
is defined to be the union of the images of
the scales $S_\Omega (A_k)$ for all $k$.
\end{defn}
\medskip

As in the proof above it follows that a group homomorphism 
$\theta : G(A_1) \to G(A_2) $
 lifts to a homomorphism in $\Omega$ if and only if 
$\theta(S_{\Omega}(A_1)) \subseteq S_{\Omega}(A_2)$.
\medskip

\begin{defn}
Let $\Omega$ be as in the last
defintion.  Then $\Omega$ is said to have the {\em uniqueness
property} if whenever $\phi, \psi : A(H) \rightarrow A$ are two
regular star-extendible 
embeddings from $\Omega$, with the same induced maps between  the dimension
distribution groups, $\tilde{\phi}, \tilde{\psi}: G(A(H)) \rightarrow
G(A)$, then $\phi$ and $ \psi$ are inner conjugate.  
\end{defn}

\medskip

We can now obtain the following abstract theorem.
\medskip

\begin{thm}
Let $\Omega$ be a family of star-extendible
regular embeddings between $H$-algebras 
and suppose that $\Omega$
has the uniqueness
property. If $\ca, \ca'$ are two systems in $\cf (\Omega)$ then $\ca$ and
$\ca '$ are regularly isomorphic if and only if there exists a dimension
distribution group isomorphism $\Phi: G(\ca) \rightarrow G(\ca')$ 
such that $\Phi(S_\Omega (\ca)) \subseteq S_\Omega (\ca ')$.
\end{thm}
\medskip

\begin{pf}
 Since $G(\ca')$ and $S_{\Omega}(\ca')$ are finitely
generated the composition
$$G(A_1)  \rightarrow G(\ca) \stackrel{\Phi}{\rightarrow} G(\ca')$$
factors as 
$$G(A_1) \stackrel{\Phi_1}{\rightarrow} G(A_{m_1}') \rightarrow G(\ca'),$$
for some $m_1$, where $\Phi_1$ is a matrix unit scale preserving abelian
group homomorphism.  
By the last lemma $\Phi_1$ lifts to a star-extendible regular
algebra homomorphism $\phi_1 $ in $\Omega$.
Similarly the composition
$$G(A'_{m_1}) \rightarrow G(\ca') \stackrel{\Phi^{-1}}{\rightarrow} G(\ca)$$
factors through $G(A_{n_1})$, for some $n_1$, and this map 
lifts to $\psi_1$.  Since $G(\psi_1 \circ \phi_1)$ agrees with 
the given group homomorphism $i$ from $G(A_1)$ 
to $G(A_{n_1})$, it follows from the
uniqueness property that we may replace $\psi_1$ by an inner conjugate
to obtain 
\[
\psi_1 \circ \phi_1 = \alpha_{n_1-1} \circ \dots \circ \alpha_2 
\circ \alpha_1.
\] 
  Continue in this way to obtain the 
desired regular isomorphism.
\end{pf}

\subsection{Classification of $T_2$ and $T_3$ systems}

We now obtain an abstract dimension distribution group classification
of the families $\cf_2$,  $\cf_3$ of regular star-extendible systems of
$T_2$-algebras and $T_3$-algebras.
This will
follow immediately from Theorem 3.3 once we demonstrate that the family
$\Omega_3$ of star-extendible regular embeddings between
$T_3$-algebras has the uniqueness property above. This is the 
assertion of Lemma 3.4.  
\medskip

\begin{thm}
The direct systems $\ca$ in $\cf_2$ and $\cf_3$ are classified up to regular
isomorphism by the scaled ordered group invariant $(G(\ca), S_{mu}(\ca))$
where $S_{mu}(\ca))$ is the matrix unit scale. 
\end{thm}

\medskip
  Observe first that a multiplicity one embedding $\theta: A_1 \rightarrow
A_2$ between $T_r$-algebras has inner conjugacy class determined by the
inner conjugacy class of the restriction $\theta|A_1 \cap A_1^*$. 
Also note that this restriction necessarily respects the ordering of the
$r$ summands of $A_1 \cap A_1^*$ and $A_2 \cap A_2^*$.  Conversely, to
any assignment of $r$ ordered objects to $r$ ordered boxes (in which
subsequent objects must not be assigned to preceding boxes) there is an
associated multiplicity one embedding $\theta$, at least if the
summands of  $A_2 \cap A_2^*$ are sufficiently large.
In the case of $T_3$-algebras it follows that
there are ten classes of mutiplicity one
embeddings, $\theta_1, \theta_2, \dots, \theta_{10}$ in which the
element $a \oplus b \oplus c$ of $A_1 \cap A_1^*$ is mapped to $A_2 \cap
A_2^*$ according to the ten ordered assignments (or partitionings)
indicated below:
\[
\theta_1: \; \; (a~ b~ c~ | 0~ 0~ 0~ | 0~ 0~ 0~)
\ \ \ \theta_2: \; \; (a~ b~ 0~ | c~ 0~ 0~ | 0~ 0~ 0~)
\]
\[
\theta_3: \; \; (a~ b~ 0~ | 0~ 0~ 0~ | c~ 0~ 0~)
\ \ \ \theta_4: \; \; (a~ 0~ 0~ | b~ c~ 0~ | 0~ 0~ 0~)
\]
\[
\theta_5: \; \; (a~ 0~ 0~ | b~ 0~ 0~ | c~ 0~ 0~)
\ \ \ \theta_6: \; \; (a~ 0~ 0~ | 0~ 0~ 0~ | 0~ b~ c~)
\]
\[
\theta_7: \; \; (0~ 0~ 0~ | a~ b~ c~ | 0~ 0~ 0~)
\ \ \ \theta_8: \; \; (0~ 0~ 0~ | a~ b~ 0~ | c~ 0~ 0~)
\]
\[
\theta_9: \; \; (0~ 0~ 0~ | a~ 0~ 0~ | b~ c~ 0~)
\ \ \ \theta_{10}: \; (0~ 0~ 0~ | 0~ 0~ 0~ | a~ b~ c~)
\]
In general one can show, by a simple recursion argument, that there are
$R = ({2r-1 \atop r})$ ordered assignments of $r$ objects to $r$ boxes.
\medskip

With the ordered $R$-tuple $\theta_1, \theta_2, \dots, \theta_R$
specified we can now define the {\em multiplicity signature} of a
regular star-extendible embedding $\phi$ between $T_r$ algebras to be
the (unique) $R$-tuple $sig({\phi}) = \{r_1, \dots, r_R\}$ for which
$\phi$ is inner conjugate to $r_1 \theta_1 \oplus \dots \oplus r_R
\theta_R$, where $r_k \theta_k$ is shorthand for $\theta_k \oplus
\dots \oplus  \theta_k$, $r_k$ times.  It is essentially a tautology
that $\phi$ that $\phi '$ are inner conjugate if and only if they have
the same signature.  It follows that {\em unital}
$T_r$-systems $\ca$ and $\ca '$ with $A_k = A_k'$ for all $k$,
are regularly isomorphic if
$sig (\alpha_k) = sig (\alpha_k')$ for all $k$.

\medskip
{\bf Lemma 3.4.} ~ {\it (Uniqueness)}\ \ \  
{\it Let $\phi, \psi: A_1 \rightarrow A_2$ be  
regular star-extendible embeddings between $T_3$-algebras.  If the
induced scaled group homomorphisms $G \phi, G \psi: G(A_1) \rightarrow
G(A_2)$ agree then $\phi, \psi$ are inner conjugate. }

\medskip
{\em Proof.} ~ The general case follows in a
straightforward way from the case $A_1 = T_3 $.  Let $x=e_{12},
y=e_{23}, z = e_{13}$ be the off-diagonal matrix units of $A_1$, and
for notational economy let $x, y, z$  also denote the corresponding elements
of $G(A_1) = T_3 (\ZZ)$.  Note that $G \theta_k : GA_1 \rightarrow
GA_2$ is determined as a $\pi$-respecting group homomorphism
by the triple $G \theta_k(x), G \theta_k(y), G
\theta_k(z)$ and, more generally, $G \phi$ is determined by the triple
$G
\phi(x), G \phi(y), G \phi(z)$.  
Each matrix $G \phi(x)$ has up to six nonzero entries,
and so the lemma will follow if we can show that the
eighteen entries for the triple $G \phi(x), G \phi(y), G \phi(z)$
determine the
multiplicity signature  $\underline{r} = 
\{r_1, r_2, \dots, r_{10}\}$
of $\phi$.  Writing these eighteen entries as a
row vector $\underline{w}$ we have 
$\underline{w} = \underline{r} X$ where $X$ is the 10 by 18 matrix
\[
\left [\begin {array}{ccccccccccccccccccccccccccc} 
                    1&0&0&0&0&0&1&0&0&0&0&0&1&0&0&0&0&0
\\\noalign{\medskip}1&0&0&0&0&0&0&1&0&0&0&0&0&1&0&0&0&0
\\\noalign{\medskip}0&1&0&0&0&0&0&0&0&1&0&0&0&1&0&0&0&0
\\\noalign{\medskip}1&0&0&0&0&0&0&0&1&0&0&0&0&0&1&0&0&0
\\\noalign{\medskip}0&1&0&0&0&0&0&0&0&0&1&0&0&0&1&0&0&0
\\\noalign{\medskip}0&0&1&0&0&0&0&0&0&0&0&1&0&0&1&0&0&0
\\\noalign{\medskip}0&0&0&1&0&0&0&0&0&1&0&0&0&0&0&1&0&0
\\\noalign{\medskip}0&0&0&1&0&0&0&0&0&0&1&0&0&0&0&0&1&0
\\\noalign{\medskip}0&0&0&0&1&0&0&0&0&0&0&1&0&0&0&0&1&0
\\\noalign{\medskip}0&0&0&0&0&1&0&0&0&0&0&1&0&0&0&0&0&1
\end {array}\right ]
\]
One can check that this matrix has has rank 10 and so 
$\underline{r}$ may be obtained from $\underline{w}$, as required.
\hfill $\Box$
\medskip

One should note that although $x$ and $y$ determine $z$
the rank distributions of $\phi(x)$ and $\phi(y)$ need not determine the
rank distribution of $\phi(z)$.  And indeed, the $10 \times 12$ matrix
resulting from overlooking the distribution of $\phi(z)$ only has rank
9.  It follows then that there are nonconjugate $\phi, \phi'$ for
which $\phi(x)$ has the same rank distribution as $\phi'(x)$ and
$\phi(y)$ has the same rank distribution as $\phi'(y)$.
\medskip

\subsection{Classifying $T_4$ systems}

Arguing as above, there are 35 conjugacy classes 
of multiplicity one algebra homomorphisms
from $T_4$ to $T_4 \otimes M_n$, for $n \geq 4$, and so the multiplicity
signature of a regular star-extendible embedding between $T_4$-algebras
is a 35-tuple.  Arguing just as in the proof of Lemma 3.4 the family
$\Omega_4$ of such embeddings has the uniqueness property if and
only if the (ordered) sextet of integral matrices
$$G \phi(e_{12}), G \phi(e_{23}), G \phi(e_{34}), G \phi(e_{13}), G
\phi(e_{24}), G \phi(e_{14})$$
determines the multiplicity signature of a regular star-extendible
 algebra homomorphism
$\phi: T_4 \rightarrow A$,
into a $T_4$-algebra $A$.  Each image $G \phi (e_{ij})$ is a $4 \times
4$ upper triangular matrix, with up to ten nonzero entries, and so the
sextet provides 60 linear equations for the 35 -tuple.  However,
(computer) calculation shows the rank of the associated $35 \times 60$
coefficient matrix to be 31 and so $\Omega_4$ fails to have the
uniqueness property.
\medskip

On the other hand, let $\Omega^-_4$ be the family of embeddings which
are regular and star-extendible and which have no multiplicity one
summand $\theta$, of degenerate type, with 
range in the self-adjoint subalgebra.  
Clearly there are four such degenerate embeddings $\theta$.
Then
(computer) calculation shows that the associated $31 \times 60$
coefficient matrix has rank 31. Thus  $\Omega_4^-$ does have the uniqueness
property and the family $\cf_4^- = \cf(\Omega_4^-)$ 
admits classification in terms of $G(\ca)$;
\medskip

\begin{thm}
 The direct systems $\ca$ of $\cf_4^-$ 
are classified by the scaled ordered dimension 
group $(G(\ca), S_{\Omega_4^-}(\ca))$.
\end{thm}

\medskip

We note that for $\ca$ in $\cf_4^-$ the group $G(\ca)$ is
of the form $\indlimit ({\ZZ}^{10}, G(\alpha_k))$
and so somewhat more computable, in principle, than $G^r_H(\ca) = 
\indlimit ( {\ZZ}^{35}, G^r_H(\alpha_k))$.

\subsection{Systems of cycle algebras}

A 4-cycle algebra is a
digraph algebra $A$ whose reduced digraph is isomorphic to the 4-cycle
digraph $D_4$.  The simplest of these is $A(D_4)$ which is the subalgebra of $M_4(\IC)$
spanned by the matrix units 
$e_{11}, e_{22}, e_{33}, e_{44}, 
e_{13},  e_{14},e_{24}, e_{23}$.
The so-called {\em rigid} embeddings between 4-cycle algebras are
those embeddings $\phi: A_1 \rightarrow A_2$ for which $\phi$ is inner
conjugate to $r_1 \theta_1 \oplus r_2 \theta_2 \oplus r_3 \theta_3
\oplus r_4\theta_4,$ where $\theta_1, \dots, \theta_4$ are multiplicity
one star-extendible embeddings whose $K_0$ matrices have the form
\[
\left[ \begin{array}{cccc}
1 & 0 & 0 & 0\\
0 & 1 & 0 & 0 \\
0 & 0 & 1 & 0\\
0 & 0 & 0 & 1 \end{array}
\right],
\left[ \begin{array}{cccc}
1 & 0 & 0 & 0\\
0 & 1 & 0 & 0\\   
0 & 0 & 0 & 1\\   
0 & 0 & 1 & 0 \end{array}
\right],
\left[ \begin{array}{cccc}
0 & 1 & 0 & 0\\
1 & 0 & 0 & 0\\
0 & 0 & 0 & 1\\
0 & 0 & 1 & 0 \end{array}
\right],
\left[ \begin{array}{cccc}
0 & 1 & 0 & 0\\
1 & 0 & 0 & 0\\
0 & 0 & 1 & 0\\
0 & 0 & 0 & 1 \end{array}
\right].
\]
These embeddings are those associated with the digraph
automorphisms of $D_4$, rather than the digraph endomorphisms.
The dimension distribution group $G(A(D_4))$ is naturally identifiable with
${\ZZ}^4 \oplus {\ZZ}^4$, and we may order the generators so that the
induced morphism
\[
G(\phi) = \left[ \begin{array}{cc}
K_0\phi & 0\\
0 & \Gamma \phi\end{array}
\right] \ : {\ZZ}^4 \to {\ZZ}^4
\]
is realised by the integral matrix
\[
\left[ \begin{array}{cccccccc}
 r_1+r_4  &  r_2+ r_3  & 0 & 0 & 0& 0 &0 &0\\
 r_2+r_3  &  r_1+r_4  & 0 & 0 &0 &0  &0 &0\\
0 & 0 &  r_1 + r_2  &  r_3 + r_4  &0 &0 & 0  &0\\
0 & 0 &  r_3 + r_4  &  r_1 + r_2  &0 &0 &0 &0 \\
0&0 &0 &0 &  r_1  &  r_4  &  r_3  &  r_2 \\
0&0 &0 &0 &  r_4  &  r_1  &  r_2  &  r_3 \\
0&0 & 0 &0 &  r_3  &  r_2  &  r_1  &  r_4 \\
0&0 &0 &0 &  r_2  &  r_3  &  r_4  &  r_1 \\
\end{array} \right] .
\]
In view of the form of $\Gamma_\phi$, in which the multiplicity 
signature for $\phi$ is evident, the family of
such embeddings has the uniqueness property.  Likewise the family
of rigid embeddings between $2n$-cycle algebras, has the uniqueness
property and it follows from Theorem 3.6 that these 
systems admit a classification
in terms of a scaled dimension distribution group.
\medskip

However, for these systems there is a more economical and revealing 
classification 
scheme which has been developed
in Power \cite{scp-book} and in Donsig and Power \cite{don-scp-2},
\cite{don-scp-3}.
Associated with a regular embedding $\phi$ between digraph algebras
there is a natural induced homology group homomorphism $H_1 \phi$,
which, in the case of 4-cycle algebras, is a group homomorphism from  ${\ZZ}$
to $ {\ZZ}$.  For a rigid embedding $\phi$, as above, this map is realised by
the $1 \times 1$ integral matrix 
$$[r_1 - r_2 + r_3 - r_4].$$
\medskip

The fact that $H_1 \phi$ and $K_0 \phi$ determine the inner conjugacy
class of $\phi$ 
is at the heart of the
classification of rigid 4-cycle systems $A$ in terms
of $K_0H_1(A)$ and various scales.  

We return to these systems in the next section.

\section{Classifications with $G^r_H(\ca)$}

Let $H$ be a connected transitive antisymmetric digraph,
so that $A(H)$ is a triangular digraph algebra with digraph $H$.
Let $G^r_H(\ca)$ be the regular Grothendieck group of the regular star-extendible
system $\ca = \{A_k, \alpha_k\}$,
with order induced by the scale $\Sigma^r_H(\ca)$. Note that if $\alpha : A \to A'$
is a regular star-extendible homomorphism of digraph algebras
then there is a natural induced scaled ordered group homomorphism
\[
\hat{\alpha} : G^r_H(A) \to  G^r_H(A')
\]
which is induced by the correspondence $V^r_H(A) \to V^r_H(A')$
given by $[\phi] \to [\alpha \circ \phi]$.
Also it follows from the definition that $G^r_H(\ca)
 = \indlimit (G^r_H(A_k), \hat{\alpha})$ so that
$G^r_H(\ca)$ is a scaled ordered dimension group.
\medskip

Let $\cf_H$ be the family of regular star-extendible systems $\{A_k, \alpha_k\}$
of $H$-algebras with the convention that $A_1 = A(H)$.
We say, simply, that $\ca, \ca'$ in $\cf_H$ are isomorphic
if there exists a commuting diagram as in section 2.1, and that 
$\ca, \ca'$ are stably isomorphic if there is such a diagram for the
stable systems 
$\{A_k \otimes M_k, \alpha_k \otimes i_k\}$,
$\{A_k' \otimes M_k, \alpha_k' \otimes i_k\}$
where $\{M_k, i_k\}$ is the system for the
stable algebra $M_\infty.$
\medskip

%
%
%
%
%
%

Recall that the inner unitary equivalence class $[\phi]$
of a regular star-extendible homomorphism
$\phi : A(H) \to A$ is determined by its multiplicity signature,
which is an $r$-tuple of nonnegative integers,
where $r$ is the number of multiplicity one embeddings from $A(H) $to $A$.
It is this basic fact that ensures that the natural inclusion
$V_H^r(A) \to G^r_H(A)$ is an injection, identifiable with the injection
${\Bbb Z}^r_+ \to {\Bbb Z}^r$,
and that $G^r_H(A)$ is a dimension group.
This basic fact also leads to the following uniqueness lemma.
\medskip

\begin{lma}
(Uniqueness.) Let $\phi, \psi : A_1 \to A_2$ be regular
star-extendible homomorphisms between
$H$-algebras with $\hat{\phi} = \hat{\psi}$.
Then $\phi$ and $\psi$ are inner conjugate.
\end{lma}
\medskip

\begin{pf}
Let $i: A(H) \to A_1$ be  a proper injection.
Since $\phi, \psi$ are inner conjugate if and only if
their restrictions to $i(A(H))$ are inner conjugate we may as well assume
that $A_1 = A(H)$.
The injection $i$ determines the element $[i]$ of $V_H^r(A_1)$,
which we may identufy with the element $[i]$
of $G^r_H(A_1)$. By definition
$\hat{\phi}([i]) = [\phi \circ i]$, the class in $G^r_H(A_1)$. Thus the hypotheses
imply that $[\phi] = [\psi]$ as classes in $G^r_H(A_2)$ and the lemma follows.
\end{pf}
\medskip

For a given ordered group homomorphism
$\gamma : G^r_H(A_1) \to G^r_H(A_2)$
there need not exist a lifting $\phi : A_1 \to A_2 \otimes M_k$
with $\hat{\phi} = \theta$. However we can guarantee this existence in
two ways.

Let $\theta_1, \dots , \theta_R$
be the distinct multiplicity one
injections $\theta_i : A(H) \to A(H) \otimes M_k$
with $k \ge R$, corresponding to the $R$ endomorphisms of the digraph
$H$.
In analogy with the definition
of the "multi-scale" $S_H(B)$ given in section 3, define the multi-cone
\[
C_H(A_1) \subseteq G^r_H(A_1)_+ \times \dots \times G^r_H(A_1)_+
\]
as the subset of the $R$-fold product
consisting of the elements
\[
([\phi \circ \theta_1], \dots ,[\phi \circ \theta_R] )
\]
associated with regular star extendible embeddings
$\phi : A(H) \to A_1 \otimes M_k$, for arbitrary $k$.
Now, if $\gamma : G^r_H(A_1) \to G^r_H(A_2)  $
preserves the multi-cone then
\[
(\gamma ([\theta_1], \dots , \gamma ([\theta_R] )) = 
([\phi \circ \theta_1]
, \dots , [\phi \circ \theta_R])
\]
for some embedding $\phi : A(H) \to A_2 \otimes M_k$.
Moreover,
since $\hat{\phi} $ agrees with $\theta$ on the 
generators of $G^r_H(A(H))$,
it follows that $ \hat{\phi} = \theta$.

Naturally, the multiscale $\tilde{\Sigma}_H(A_1)$ is
defined as the subset of the
multi-cone associated with the classes of embeddings of $A(H)$ into
the algebra $A_1$ itself.

An alternative equivalent requirement that guarantees the existence of a
lifting for $\gamma$ is to take into account
the fact that an induced
ordered group homomorphism $\hat{\phi}$ respects 
the natural right action
of the semigroup $End(H) = \{\theta_1, \dots , \theta_k\}$
and to require that $\gamma$ respects this right action.
We discuss this important perspective in \cite{scp-afa}.

\begin{lma} (Existence.) Let $A_1, A_2$ be $H$-algebras
and let $\theta : G^r_H(A_1) \to G^r_H(A_2)$  be an ordered group homomorphism
which preserves the multicone.
Then for sufficiently large $k$ there exists a regular star-extendible
algebra homomorpism
$\phi : A_1 \to A_2 \otimes M_k$ such that $\hat{\phi} = \theta $.
\end{lma}
\medskip

In view of the lemmas above we can obtain the following theorem
in the usual fashion.

\begin{thm}
Systems $\ca, \ca'$ in $\cf_H$
are stably regularly isomorphic if and only if the scaled ordered
dimension groups $G^r_H(\ca) $ and $G^r_H(\ca')$ are isomorphic
by an isomorphism which preserves the multiscale.
\end{thm}
\medskip

%
%
%
%

\noindent {\bf Examples.}
\medskip

For $H = T_3$ the regular Grothendieck group invariant, which has
the form $\indlimit ({\Bbb Z}^{10}, \hat{\alpha_k})$ , and
the integral matrix for
$\alpha_k$ can be computed from the multiplicity
signature and the $10 \times 10$ multiplication
table for the semigroup of order preserving partitions.
\medskip

For subfamilies of $T_3$-algebra systems it becomes 
 appropriate to consider a reduced
Grothendieck group invariant as in the following example.
\medskip

Let $ \cg$ be the family of embeddings $\phi$ between $T_3$-algebras which
are unitarily equivalent to a direct sum
$\rho_r \oplus \sigma_{3s} \oplus \delta_t$
where $\rho_r$ is a refinement embedding of multiplicity
$r$,
\medskip

\[
\rho_r : 
\left[ \begin{array}{ccc}
a&e&g\\
 &b&f\\
 & &c
\end{array} \right]
\to 
\left[ \begin{array}{ccc}  
a \otimes I_r&e\otimes I_r&g\otimes I_r\\    
 &b\otimes I_r&f\otimes I_r\\
 & &c\otimes I_r
\end{array} \right] 
\]

\medskip

\noindent where $\sigma_{3s}$ is a standard embedding of multiplicity $3s$,

\medskip

\[
\sigma_{3s} : [a] \to
\left[ \begin{array}{ccc}
\sigma_s(a)&0&0\\
 &\sigma_s(a)&0\\ 
 & &\sigma_s(a)
\end{array} \right],
\]

\medskip

\noindent and where $\delta_t$ is a degenerate embedding of multiplicity $t$,
with
\medskip

\[
\delta_t : a \to
\left[ \begin{array}{ccc}
\sigma_t(a)&0&0\\
 &0&0\\
 & &0  
\end{array} \right].
\]

\medskip

View the embedding $\phi$ as having reduced multiplicity signature
$\{r, s, t\}$
and view $\phi$ itself as a direct sum of the basic embeddings
$\rho_1, \sigma_3$ and $ \delta_1$.
The family $\cg$ is closed under compositions and the
corresponding reduced Grothendieck group $G_{\cg}(\ca)$ for the system
$\ca = \{A_k, \alpha_k\}$, with $sig(\alpha_k) = \{r_k, s_k, t_k\}$
may be identified with the dimension group

\medskip
\[
\indlimit ({\Bbb Z}^3, \left[ \begin{array}{ccc}
r_k&0&0\\
s_k &r_k+s_k&0\\
t_k &t_k &r_k+s_k+t_k  
\end{array} \right]
).
\]

\medskip

These dimension groups and their multiscales
classify the stable systems in the family up to regular isomorphism.
\medskip

Let us also illustrate these methods with an application to
the locally finite
algebras $A$ of the towers $\ca = \{A_k, \alpha_k\} $
of $2m$-cycle algebras where the inclusions
$\alpha_k : A_k \to A_{k+1}$ are of rigid
type, as indicated in section 3.4. Let $V_m(\ca)$ be the
semigroup in $V_{D_{2m}}^r(A)$ corresponding to classes $[\phi]$ associated
with rigid embeddings and let $G_m(\ca)$ be the associated 
Grothendieck group.
Plainly $G_m(A_k)$ is identifiable with ${\Bbb Z}^{2m}$ with generators
corresponding to the $2m$ symmetries of the digraph $D_{2m}$.
Furthermore, $G_m(\alpha_k)$, the induced ordered group homomorphism,
is effected by viewing ${\Bbb Z}^{2m}$ as the group ring, ${\Bbb Z}[D_{2m}]$
say,
for the automorphism group of $D_{2m}$ and by regarding $G_m(\alpha_k)$
as left multiplication by $g(\alpha_k) \in {\Bbb Z}[D_{2m}]$, the
element $r_1\theta_1 + \dots + r_{2m}\theta_{2m}$ of ${\Bbb Z}[D_{2m}]$
for the multiplicity signature of $\alpha_k$.
In case $m = 3$ this is given by multiplication by the matrix
\medskip

\[
\left[ \begin{array}{cccccc}
r_1&r_2&r_5&r_4&r_3&r_6\\
r_2&r_1&r_4&r_3&r_6&r_5\\
r_3&r_6&r_5&r_2&r_1&r_4\\
r_4&r_5&r_6&r_1&r_2&r_3\\
r_5&r_4&r_1&r_6&r_3&r_2\\
r_6&r_3&r_2&r_5&r_4&r_1
\end{array} \right].
\]

\medskip

In this way one obtains an identification of the reduced Grothendieck group
$G_m(\ca)$ as

\medskip

\[
G_m(\ca) = \indlimit({\Bbb Z}[D_{2m}], g(\alpha_k)).
\]

\medskip

We remark that one can see from Section 3.4 in the  case of 4-cycle algebras,
that the group $G_m(\ca)$ appears naturally as a summand of the
dimension distribution group $G(\ca)$.
\medskip

A key result in  Donsig and Power \cite{don-scp-3} is that for $m \ge 3$ the
locally finite algebras
$A= \indlimit (A_k, \alpha_k)$, $A' = \indlimit (A_k, \alpha_k)$
of such towers $\ca, \ca'$ are star-extendibly isomorphic if and only if the
towers are regularly isomorphic. In view of this we may define $G_m(A)$ in a
well-defined way as the group $G_m(\ca)$. Combining
these fact with Theorem 4.3 (adapted to the reduced Grothendieck group)
leads to the following classification. This  classification may 
also be extended to the operator algebras of the systems.
\medskip

\begin{thm}
Let $m \ge 3$ and let $A$ and $A'$ be the locally finite complex
algebras associated with the towers $\ca = \{A_k, \alpha_k\},$ 
$\ca' = \{A_k', \alpha_k'\}$
consisting of $2m$-cycle algebras, for $m \ge 3,$
and rigid embeddings. Then $A$ and $A'$ are
stably isomorphic if and only if the ordered groups
$G_m(A), G_m(A')$ are isomorphic by an isomorphism
which preserves the multicone.
\end{thm}

\medskip
 
Although we have focused on 
systems $\{A_k, \alpha_k\}$
for which the reduced digraph of $A_k$ is connected, with bounded diameter,
the arguments and results of the paper
also extend in a routine way to more general systems, and in particular,
to systems of direct sums of $H$-algebras.

\end{document}